\documentclass[a4paper,12pt]{amsart}

\usepackage[latin1]{inputenc}
\usepackage[T1]{fontenc}
\usepackage[english]{babel}

\usepackage{mathrsfs}
\usepackage{amscd}
\usepackage{amsfonts}
\usepackage{amsmath}
\usepackage{amssymb}
\usepackage{amstext}
\usepackage{amsthm}
\usepackage{amsbsy}

\usepackage{xspace}
\usepackage[all]{xy}
\usepackage{graphicx}
\usepackage{url}
\usepackage{latexsym}

\theoremstyle{definition}

\title[]{Murray A. Marshall [24.3.1940 - 1.05.2015]\\  A brief mathematical obituary}
\author{by Salma Kuhlmann}
\address{Fachbereich Mathematik und Statistik,\newline \indent
Universit\"at Konstanz,\newline \indent
78457 Konstanz, Germany.}
\email{salma.kuhlmann@uni-konstanz.de}
\date{\today}
\begin{document}
\maketitle
\section{Introduction}
 On April 22 2015,  Murray Marshall was invited by Igor Klep and Victor Vinnikov to talk in the special session on ``Operator Theory, Real Algebraic Geometry, and Moment Problems'' that they co-organised for IWOTA 2015 in  Tbilisi, Georgia, July 6--10. This invitation to an IWOTA edition was considered by the organisers as long overdue, given Marshall's important contributions to Real Algebraic Geometry in general, and to the multi-dimensional Moment Problem in particular. Marshall declined on April 27th writing:

\smallskip

{\it Dear Victor,
Thanks very much for the invitation. I am sorry but I must decline. The
reason is that I don't enjoy traveling as much as I used to. North
America or Europe is about the maximum distance I am willing to deal
with. I hope you are well. Murray}

\smallskip

This response startled us, as Marshall has always been a grand traveler when it came to mathematics. Our fears were alas confirmed;
Marshall passed away just a few days later. In July, I held a talk at that session myself on joint work with Marshall, and was asked by the editors of the proceedings volume to write this mathematical obituary.

\smallskip

Marshall has been updating  his cv (see below) and his homepage [http://math.usask.ca/~marshall/index.html]
 regularly. These contain links to the pdf files of his publications, as well as a listing of supervised graduate students.
According to these valuable sources of information, Murray classified his own research contributions into two main areas:
1. Positive Polynomials and Sums of Squares and 2. Orderings, Valuations and Quadratic Forms. In this brief synopsis,
we will focus on our joint contributions to the multi-dimensional and infinite-dimensional Moment Problem (which mainly belongs to his first main research area).

\section{On the moment problem}

In the winter of 1999-2000 Murray gave a series of lectures in our Algebra Seminar at the University of Saskatchewan.
We cite the following from the Preface of his little blue book:

\smallskip

{\it The goal of these lectures was to provide an elementary introduction (at a level suitable for first year graduate students)
to recent developments in real algebraic geometry related to Schm\"udgen's solution of the Moment Problem in 1991.}

\smallskip

The little blue book (Positive Polynomials and Sums of Squares, Dipartimento di Matematica Dell' Universita Di Pisa,
Dottorato di Ricerca In Matematica 2000) not only served for training our graduate students, but also provided later
the core material for his - by now famous and widely used - monograph {\it Positive polynomials and sums of squares} in the AMS Math. Surveys and Monographs series.

\smallskip

After that memorable seminar, we decided to surprise Murray by organizing a Colloquiumfest at the University of Saskatchewan in the honour
of his 60th birthday [24. 03. 2000] and invited Konrad Schm\"udgen as a key note speaker. Schm\"udgen stated two conjectures:
the solvability of the Moment Problem for cylinders with compact base in any dimension, and the saturation of the preorder
associated to the two dimensional strip\footnote{Consider the two-dimensional strip $K_S$, described as a semi-algebraic set in the real plane by  $S: =\{x, 1-x\}$. Show that every polynomial $p(x, y)$ which is non-negative on $K_S$ belongs to the associated quadratic preorder $T_S$.}. It was in the summer of 2000 in Saskatoon that, challenged by Schm\"udgen's conjectures,
our long and fruitful collaboration was born. By the end of the summer we had solved the first conjecture.
The second was established by Murray alone several years later.\footnote{{\it Polynomials non-negative on a strip}, Proceedings AMS, 138 (2010) 1559--1567. The referee of this Obituary wrote in his report: I remember that this strip theorem was quite a surprise and made a big splash;
Schm\"udgen and others have seriously pursued this but in the end it was Murray's talent that succeeded.}

\smallskip

Murray's interest in the Moment Problem only grew thereafter. While supervising  MSc student W. Fan in 2004
we discovered N. I. Akhiezer's book on The Classical Moment Problem. [Murray cherished this book so much, that he
re-assigned readings from it to our current MSc Student A. I. Jinadu, more than ten years later].
Murray thus discovered very early A. E. Nussbaum's and T. Carleman's contributions to the Moment Problem.
He seemed fascinated by the challenge of circumventing analytic and operator theoretic proofs and providing
instead algebraic ones (based on his ideas using appropriate localizations of polynomial algebras).
This is evident in his most recent publications in  Math. Scandinavica (part II is to appear posthumous).

\smallskip

In the last five years, Murray became excited about the relationship between the topological closure of quadratic modules
with respect to various locally convex topologies on the one hand, and the solvability of the Moment Problem for continuous
linear functionals on the other. Our intense collaboration together with Mehdi Ghasemi - one of Murray's favourite PhD students -
culminated in the discovery of applications to the infinite dimensional Moment Problem, and in particular to symmetric
algebras of locally convex topological vector spaces (submitted, together with Maria Infusino).
In this context, Murray became increasingly intrigued by the question of determinacy of the solving measure,
as witnessed by the last note that he sent us about this issue (submitted to this proceedings volume).

\smallskip

There are still many unsolved open questions related to our most recent collaboration with Murray.\footnote{For a detailed account of these see {\it Infinite dimensional moment problem: open questions and applications} by M. Infusino and S. Kuhlmann, arXiv::1610.01659.}
We will continue "digging" until we understand. It is appropriate to end by citing from an e-mail message
that Murray sent us on April 10 2015 in discussing referee comments on our paper
 (I removed references for confidentiality):

\smallskip

{\it Finally, I should say that I have no idea how the result in [x]  is
related to the [y] result in our other paper (although this, of
course, is another matter). It seems though that the result in [x] is for general nuclear spaces, whereas the result of [y]  is for
nuclear spaces of a very special sort. I hope eventually to understand all these things. Murray}
\section
{Applications to polynomial optimisation}
Application of the Moment Problem to Polynomial Optimisation is one of the
 most beautiful and useful achievements of the area
 which has been developed extensively in past two decades. Murray has worked
 on this problem and published several articles.
 Since 2012, together with M. Ghasemi, he introduced a new approach and published a
 series of articles based on geometric programming
 to provide a lower bound for the global minimum of a polynomial where no
 other methods are able to do so. His latest work in this
 series, including a collaboration with J. B. Lasserre,  generalises the geometric programming method to constrained Polynomial
 Optimisation problems and the last one, which deals with optimisation over
 semi-algebraic sets is still under review.
\section
{Short Curriculum Vitae of Murray Marshall}
\noindent\textbf{Name:} Murray Angus Marshall\\
\noindent\textbf{Date and Place of Birth:} March 24, 1940 -- Hudson Bay, Saskatchewan.\\
\noindent\textbf{Date and Place of Death:} May 1st-- Saskatoon, Saskatchewan.\\
\subsection{Main research interests}\vspace{-0.1cm}
Quadratic Forms, Spaces of Orderings, Real Algebraic Geometry, Positive Polynomials and Sums of Squares, Moment Problems, Polynomial Optimisation.
\subsection{Academic credentials}
\begin{itemize}
\item1966--\emph{M.A.\! in Mathematics, Algebra} at University of Saskatchewan.
\item1969--\emph{Ph.D.\! in Mathematics, Algebra} at Queen's University, Kingston.
\end{itemize}
\subsection{Appointments and promotions}
Department of Mathematics, University of Saskatchewan, Canada:
\begin{itemize}
\item Assistant Professor, 1969 - 1973.
\item Associate Professor, 1973 - 1978.
\item Full Professor, 1978 - 2006.
 \item Professor Emeritus, 2006 - 2015.
 \end{itemize}
\subsection{Recent research grants}
Natural Sciences and Engineering Research Council of Canada Discovery Grant, renewed two times since 2006, once in 2008 and once in 2013, renewed for 5 years in each case, for the amount of \$15,000 per year in each case.
\subsection{Graduate Students since 2006:}
\begin{itemize}
\item  Ayoola Isaac Jinadu (MSc, 2014 - 2015)
\item Mehdi Ghasemi (PhD, 2009 - 2012)
\item Pawe\l\ G\l adki (PhD, 2003 - 2007)
\item Wei Fan (MSc, 2004 - 2006)
 \end{itemize}
\subsection{Postdoctoral Fellows since 2006:}
\begin{itemize}
\item Pawe\l\ G\l adki, Sept.- Dec. 2014
\item Mehdi Ghasemi, June - Sept. 2014
\item Sven Wagner, 2010- 2011
\item Pawe\l\  G\l adki, July - Sept. 2010
\item Katarzyna Osiak, Aug.- Oct. 2008
\item Tim Netzer, May - Oct. 2007
\item Andreas Fischer, 2006- 2008
\item Igor Klep, May - Sept. 2006
 \end{itemize}
\section
{ Murray Marshall's selected publications}
\subsection{Positive Polynomials and Moment Problems}
\begin{enumerate}
\item M. Ghasemi, M. Infusino, S. Kuhlmann, M. Marshall, \emph{Representation of a continuous linear functional on subspaces of $\mathbb{R}$-algebras}, in preparation.
\item M. Ghasemi, M. Marshall,  \textit{Lower bounds for a polynomial on a basic closed semialgebraic set using geometric programming}, arXiv:1311.3726, submitted.
\item M. Ghasemi, M. Infusino, S. Kuhlmann, M. Marshall, \emph{Moment problem for symmetric algebras of locally convex spaces}, arXiv:1507.06781, submitted.
\item M. Infusino, S. Kuhlmann, M. Marshall, \emph{On the determinacy of the moment problem for symmetric algebras of a locally convex space}, arXiv:1603.07747, submitted.
\item M. Marshall, \textit{Application of localization to the multivariate moment problem II},  arXiv:1410.4609, to appear in Math. Scandinavica.
\item M. Ghasemi, S. Kuhlmann, M. Marshall, \emph{Moment problem in infinitely many variables}, Israel Journal of Mathematics, 212 (2016) 989-1012.
\item M. Ghasemi, S. Kuhlmann, M. Marshall, \textit{Application of Jacobi's representation theorem to locally multiplicatively convex topological R-algebras}, J. Functional Analysis, 266 (2014), no. 2, 1041--1049.
\item M. Ghasemi, J.B. Lasserre, M. Marshall, \textit{Lower bounds on the global minimum of a polynomial}, Computational Optimization and Applications, 57 (2014) 387--402.
\item M. Marshall, \textit{Application of localization to the multivariate moment problem}, Math. Scandinavica, 115 no. 2 (2014) 269--286.
\item M. Ghasemi, M. Marshall, Sven Wagner, \emph{Closure of the cone of sums of $2d$-powers in certain weighted $\ell_1$-seminorm topologies}, Canad. Math. Bull., 57, no 2, (2014) 289-302.
\item M. Ghasemi, M. Marshall, \textit{Lower bounds for polynomials using geometric programming}, SIAM Journal on Optimization, 22 (2012) 460--473.
\item M. Marshall, T. Netzer, \textit{Positivstellens\"atze for real function algebras}, Math. Zeitschrift, 270 (2012) 889--901.
\item J. Cimpri\v c, M. Marshall, T. Netzer, \textit{Closures of quadratic modules}, Israel J. Math., 189 (2011) 445--474.
\item J. Cimpri\v c, M. Marshall, T. Netzer, \textit{On the real multidimensional rational K-moment problem}, Transactions AMS, 363 (2011) 5773--5788.
\item M. Ghasemi, M. Marshall, \textit{Lower bounds for a polynomial in terms of its coefficients}, Archiv der Mathematik, 95 (2010) 343--353.
\item M. Marshall, \textit{Polynomials non-negative on a strip}, Proceedings AMS, 138 (2010) 1559--1567.
\item J. Cimpri\v c, S. Kuhlmann, M. Marshall, \textit{Positivity in power series rings}, Advances in Geometry, 10 (2010) 135--143.
\item M. Marshall, \textit{Representation of non-negative polynomials, degree bounds and applications to optimization}, Canad. J. Math., 61 (2009) 205--221.
\item M. Marshall, \textit{Positive polynomials and sums of squares}, AMS Math. Surveys and Monographs 146 (2008) 187+xii pages.
\item M. Marshall, \textit{Representation of non-negative polynomials with finitely many zeros}, Annales de la Faculte des Sciences Toulouse, 15 (2006) 599-609.
\item M. Marshall, \textit{Error estimates in the optimization of degree two polynomials on a discrete hypercube}, SIAM Journal on Optimization,16 (2005)  297-309.
\item S. Kuhlmann, M. Marshall, N. Schwartz, \textit{ Positivity, sums of squares and the multi-dimensional moment problem II}, Advances in Geometry, 5 (2005) 583-606.
\item M. Marshall, \textit{Approximating Positive Polynomials Using Sums Of Squares}, Can. math. bulletin, 46 (2003) 400-418.
\item M. Marshall, \textit{Optimization of Polynomial Functions}, Can. math. bulletin, 46 (2003) 575-587.
\item S. Kuhlmann, M. Marshall, \textit{Positivity, sums of squares and the multi-dimensional moment problem}, Trans. Amer. Math. Soc.  354 (2002), 4285-4301.
\item M. Marshall, \textit{A General Representation Theorem For Partially Ordered Commutative Rings},  Math. Zeitschrift 242 (2002), 217-225.
\item M. Marshall, \textit{Extending The Archimedean Positivstellensatz To The Non-Compact Case}, Can. math. bulletin, 44 (2001) 223-230.
\item M. Marshall, \textit{A Real Holomorphy Ring Without The Schm\"udgen Property}, Can. math. bulletin, 42 (1999) 354-35.
\end{enumerate}
\subsection{Orderings, Valuations and Quadratic Forms}
\begin{enumerate}
\item P. Gladki, M. Marshall, \textit{Witt equivalence of function fields of curves over local fields},  arXiv: 1601.08085.
\item P. Gladki, M. Marshall, \textit{Witt equivalence of function fields over global fields}, to appear in Trans. Amer. Math. Soc.
\item P. Gladki, M. Marshall, \textit{Quotients of index two and general quotients in a space of orderings}, Oberwolfach Preprint 2011-36, Fundamenta Mathematicae, 229 (2015) 255-275.
\item P. Gladki, M. Marshall, \textit{Orderings and signatures of higher level on multirings and hyperfields}, J. K-Theory, 10 (2012) 489-518.
\item S. Kuhlmann, M. Marshall, K. Osiak, \textit{Cyclic 2-structures and spaces of orderings of power series fields in two variables}, J. Algebra, 335 (2011) 36-48.
\item M. Machura, M. Marshall, K. Osiak, \textit{Metrizability of the space of R-places of a real function field}, Math. Zeitschrift, 266 (2010) 237-242.
\item D. Gondard, M. Marshall, \textit{Real holomorphy rings and the complete real spectrum}, Annales de la Faculte des Sciences Toulouse 19 (2010), Fascicule Special, 57-74.
\item P. Gladki, M. Marshall, \textit{On families of testing formulae for a pp formula}, Cont. Math. 493 (2009) 181-188.
\item P. Gladki, M. Marshall, \textit{The pp conjecture for the space of orderings of the field $R(x,y)$}, J. Pure \& Applied Algebra, 212 (2008) 197-203.
\item P. Gladki, M. Marshall, \textit{The pp conjecture for spaces of orderings of rational conics}, Algebra and its Applications, 6 (2007) 245-257.
\item J. Cimpric, M. Kochetov, M. Marshall, \textit{Orderability of pointed cocommutative Hopf algebras}, Algebras and Representation Theory (2007) 25-54.
\item M. Marshall, \textit{Real reduced multirings and multifields}, J. Pure  \& Applied Algebra, 205 (2006) 452-468.
\item M. Marshall, \textit{Local-global properties of positive primitive formulas in the theory of spaces of orderings}, J. Symbolic Logic, 71 (2006) 1097-1107.
\item M. Dickmann, M. Marshall, F. Miraglia, \textit{Lattice ordered reduced special groups}, Annals of Pure and Applied Logic, 132 (2005) 27-49.
\end{enumerate}
\subsection{Other Topics}
A. Fischer, M. Marshall, \textit{Extending piecewise polynomial functions in two variables}, Annales de la Faculte des Sciences Toulouse, 22 (2013) 253-268.
\section{Acknowledgements}  I wish to thank Mehdi Ghasemi for providing the paragraph on optimisation and Maria Infusino for editing the list of publications.

\end{document}